\documentclass{amsart}
\usepackage{amssymb}
\theoremstyle{definition}
\newtheorem{thm}{Theorem}
\newtheorem{lem}{Lemma}
\newtheorem{cor}{Corollary}
\begin{document}
\title[Sums of like powers and some dense sets]{SUMS OF LIKE
POWERS AND SOME DENSE SETS}
\author{\v Zarko Mijajlovi\' c, Milo\v s Milo\v sevi\' c,
Aleksandar Perovi\' c}
\address{Facylty of Mathematics, University of Belgrade,
Studentski trg 16, 11000 Belgrade, Serbia and Montenegro}
\email{zarkom@eunet.yu, mionamil@eunet.yu, pera@sf.bg.ac.yu}
 \maketitle
\begin{abstract}
In this paper we introduce the notion of the $P$-sequences and apply their properties
in studying representability of real numbers. Another application of $P$-sequences
we find in generating the Prouhet-Tarry-Escott pairs.

\end{abstract}
\section{Introduction}

First we define the notion of $P$-sequence and establish their
basic properties which will be used in later sections. In the next section  we
use the properties of $P$-sequences in the study of representability of real numbers by
sequences of reals. Recall that a real number $r$ is representable by a sequence
$\langle a_n\ |\ n\in\mathbb N\rangle$ if there is $S\subseteq \mathbb N$
such that
$
r=\sum_{n\in S}a_n.
$

The following  result on representability of real numbers is due to Kakeya:
Suppose that $A=\langle a_n\ |\ n\in\mathbb N\rangle$ is a decreasing sequence
of positive reals which converges to $0$ and that $s=\sum a_n$, $0< s\leqslant+\infty$.
Then the following assertions are equivalent:
\begin{itemize}
\item Each $r\in (0,s]$ is representable by means of $A$;
\item $a_n\leqslant\sum_{k=n+1}^\infty a_k$, for each $n$.
\end{itemize}
Though we do not use Kakeya's theorem in our proofs, some particular cases of our examples are
its consequences. However, the most interesting cases cannot be obtained
by it.

Finally, using the $P$-sequences, we obtain  new methods of
generating sums of like powers i.e. the Prouhet-Tarry-Escott pairs.

\section{$P$-sequences}

We use symbols $\mathbb{Z}$, $\mathbb{N}$, $\mathbb{N}^+$,
$\mathbb{R}$ and $\mathbb{R}^+$ to represent the sets of integers,
nonnegative integers, positive integers, real numbers and positive
real numbers, respectively. In addition, we also adopt the
convention that $0^0=1$.
\vspace{1ex}

The notion of a $P$-sequence is recursively defined as follows:
\begin{itemize}
\item $\langle 1,-1\rangle$ is a $P$-sequence; \item If $\langle
a_0,\dots, a_k\rangle$ is a $P$-sequence and $a_0=-a_k$, then
$$
\langle a_0,\dots,a_k,a_k,\dots,a_0\rangle
$$
is also a $P$-sequence; \item If $\langle a_0,\dots, a_k\rangle$
is a $P$-sequence and $a_0=a_k$, then
$$
\langle a_0,\dots,a_{k-1},0,-a_{k-1},\dots,-a_0\rangle
$$
is also a $P$-sequence;
\item Each $P$-sequence can be obtained
only by finite use of the above clauses.
\end{itemize}
We denote the $n$-th $P$-sequence by $P_n$ (assuming that they are
ordered by their increasing lengths). For instance,
$$
P_1=\langle 1,-1\rangle,\  P_2=\langle 1,-1,-1,1\rangle,\
P_3=\langle 1,-1,-1,0,1,1,-1\rangle\ {\rm etc}.
$$
For an arbitrary positive integer $n$, the $n$-th $P$-sequence $P_n
=\langle a_0,\dots,a_k\rangle$ and any integer $s\geqslant 0$ let
us define a polynomial function $F_{n,s}(x)$ over $\mathbb{R}$ by
$$
F_{n,s}(x)=\sum_{i=0}^ka_i(i+x)^s.
$$
\begin{lem}
\label{1}
$F_{n,s}\equiv 0$ for $s=0,\dots,n-1$.
\end{lem}
\noindent{\bf Proof.} Since
$$
F_{n,s}(x)=((n-1)\cdots (n-1-s))^{-1}F_{n,n-1}^{(n-1-s)}(x), \qquad 0\leq s < n-1,
$$
where $F_{n,n-1}^{(n-1-s)}$ is the $(n-1-s)$-th derivative of
$F_{n,n-1}$, it is sufficient to prove that
\begin{eqnarray}
\label{2}
F_{n,n-1}\equiv 0.
\end{eqnarray}
We prove the lemma by induction on $n$. Trivially (1) is true for
$n=1$, so let us assume that for some $n\geqslant 1$ the equality (1)
holds . We have the following two cases:
\begin{itemize}
\item $n=2m$. Assuming that $P_{2m}=\langle a_0,\dots,a_k\rangle$,
we have that
$$
P_{2m+1}=\langle a_0,\dots,a_{k-1},0,-a_{k-1},\dots,-a_0\rangle
$$
and
\begin{eqnarray*}
F_{2m+1,2m}(x)&=&\sum_{i=0}^{k-1}a_i(i+x)^{2m}-\sum_{i=0}^{k-1}a_i(2k-i+x)^{2m}\\
&=&\sum_{i=0}^ka_i(i+x)^{2m}- \sum_{i=0}^ka_i(i-2k-x)^{2m}\\
&=& F_{2m,2m}(x)-F_{2m,2m}(-2k-x).
\end{eqnarray*}
Then
\begin{eqnarray*}
F_{2m+1,2m}'(x)&=& 2m\underbrace{F_{2m,2m-1}(x)}_{=0}
 + 2m\underbrace{F_{2m,2m-1}(-2k-x)}_{=0} \\
&=&0,
\end{eqnarray*}
so $F_{2m+1,2m}$ is a constant function. Since
$F_{2m+1,2m}(-k)=0$, we conclude that $F_{2m+1,2m}\equiv 0$.
\item $n=2m+1$. Similarly to the previous case one can easily check that
$F_{2m+2,2m+1}$ is a constant function. Since
$F_{2m+2,2m+1}(-\frac{2k+1}{2})=0$, we conclude that
$F_{2m+2,2m+1}\equiv 0$ as well. \hfill $\square$
\end{itemize}
\begin{thm}
\label{3}
For $s\geqslant n$ the degree of $F_{n,s}$ is equal to $s-n$.
\end{thm}
{\bf Proof.} Clearly, it is sufficient to prove that
\begin{eqnarray}
\label{4}
F_{n,n}\equiv const.\not=0.
\end{eqnarray}
Observe that an immediate consequence of lemma \ref{1} is the fact that
each $F_{n,n}$ is a constant function.

The proof goes by induction on $n$. $F_{1,1}\equiv -1$, so let us
assume that for some $n\geqslant 1$ the relation (\ref{4}) holds.
\begin{itemize}
\item $n=2m$. Assuming that $P_{2m}=\langle a_0,\dots,
a_k\rangle$, we have
$$
P_{2m+1}=\langle a_0,\dots,a_{k-1},0,-a_{k-1},\dots,-a_0\rangle
$$
and
\begin{eqnarray*}
F_{2m+1,2m+1}(x)&=&
\sum_{i=0}^{k-1}a_i(i+x)^{2m+1}-\sum_{i=0}^{k-1}
a_i(2k-i+x)^{2m+1}\\
&=& F_{2m,2m+1}(x)+F_{2m,2m+1}(-2k-x).
\end{eqnarray*}
$F_{2m+1,2m+1}$ is a constant function, so
$$
F_{2m+1,2m+1}(x)=F_{2m+1,2m+1}(-k)=2F_{2m,2m+1}(-k).
$$
Since $a_i=a_{k-i}$, we have
$$
F_{2m,2m+1}(-k)=-F_{2m,2m+1}(0).
$$
By the induction hypothesis $F_{2m,2m+1}$ is a linear function thus
1--1, hence
$$
F_{2m,2m+1}(-k)\not=0.
$$
\item $n=2m+1$. Similarly to
the previous case one can easily deduce that
$$
F_{2m+2,2m+2}(x)=2F_{2m+1,2m+2}(-(2k+1)/2)
$$
and
$$
F_{2m+1,2m+2}(-(2k+1)/2)=-F_{2m+1,2m+2}(1/2),
$$
which combined with the induction hypothesis implies that
$$
F_{2m+1,2m+2}(-(2k+1)/2) \not=0.
$$
\hfill $\square$
\end{itemize}

\begin{cor}
\label{5}
Let $P_n=\langle a_0,\dots,a_k\rangle$ be a $P$-sequence. Then:
\begin{enumerate}
\item $\sum\limits_{i=1}^ka_ii^s=0,\ s=0,\dots,n-1$;
\item $\operatorname{sgn}\sum\limits_{i=1}^ka_ii^n=(-1)^n$.
\end{enumerate}
\end{cor}
\noindent{\bf Proof.} (1) is an immediate consequence of theorem \ref{3}, while
(2) can be obtained by slight modification of the corresponding proof
of theorem \ref{3}.
\hfill $\square$

\section{Dense-expandable sequences}

Through this section $E$ will denote some denumerable sequence of positive real
numbers, $E(n)$ its $n$-th member, $\sum E=_{\rm def}\sum\limits_{n=0}^
\infty E(n)$, $\lim E=_{\rm def}\lim\limits_{n\to\infty}E(n)$,
$r+sE=_{\rm def}\langle r+sE(n)\ |\ n\in\mathbb{N}\rangle$ and $\mu_E$ is
a measure on $\mathbb{N}$ defined by
$$
\mu_E(S)=\sum_{n\in S}E(n),\ \ S\subseteq \mathbb{N}.
$$
We say that measure $\mu_E$ is continuous if for each $r\in[0,\infty]$
there is $S\subseteq \mathbb{N}$ such that $\mu_E(S)=r$.
\vspace{2ex}

Sequence $E$ is \emph{dense-expandable} if the set
\begin{eqnarray*}
\label{6}
\operatorname{X}(E)=_{\rm def}\{\sum_{n\in S}\varepsilon_nE(n)\ |\
\mbox{$S$ ranges over finite subsets of $\mathbb N$ and $\varepsilon_n\in\{-1,1\}$}\}
\end{eqnarray*}
is dense in $\mathbb{R}$.
\vspace{2ex}

If $\sum E$ is convergent, then
$\operatorname{X}(E)$ is bounded in $\mathbb{R}$, so it cannot be
dense. On the other hand, if $\mu_E$ is continuous, then $E$ is
obviously dense-expandable.
\begin{thm}
\label{8}
Suppose that $E$ is a sequence of positive real numbers such that $\lim E=0$
and $\sum E=\infty$. Then $\mu_E$ is continuous.
\end{thm}
\noindent{\bf Proof.}
For the fixed positive real number $c$ let $\mathcal C$ be the family of all subsets
$S$ of $\mathbb{N}$ such that $\mu_E(S)\leqslant c$. Note that $\mathcal C\not=\emptyset$
since $\lim E=0$. $\langle \mathcal C,\subseteq\rangle$ is clearly a poset, so it
has a maximal chain, say $M$. It is easy to see that
$\bigcup M\in\mathcal C$. Suppose that $\mu_E(\bigcup M)=b<c$.
Since $\mu_E(\mathbb{N})=\infty$, $\mathbb{N}\setminus \bigcup M$ is infinite.
Now $\lim E=0$ implies that there is an index $n\in \mathbb{N}\setminus \bigcup M$
such that $b+E(n)<c$. But this implies that $\mu_E(\bigcup M)<\mu_E(
\{E(n)\}\cup\bigcup M)<c$, which contradicts the maximality of $M$.
\hfill $\square$
\vspace{3ex}

A converse implication need not be true. Namely, the sequence
$$
E=\langle 1,2,2^{-1}, 2^2,2^{-2}, 2^3,2^{-3},\dots\rangle
$$
clearly generates a continuous measure $\mu_E$, but it is not convergent.
However, its limes inferior $\underline{\lim} E$ is equal to 0. Note also that
from $\underline{\lim} E=0$ and $\sum E=\infty$ does not follow necessarily
the continuity of $\mu_E$.

\begin{thm}
 Let $\varphi(n)$ be the Euler function. Then
   $A= \{\varphi(n)/n \colon n\in N\}$ is dense in
   the real interval $[0,1]$.
\end{thm}
\noindent{\bf Proof.}
First,  we remind the reader that
   $\varphi(n)/n= (1-1/p_1)\ldots (1-1/p_k)$, where
   $p_1,\ldots, p_k$ are all prime factors of $n$.
   Now, let $a_n= -\log(1-1/p_n)$, where
   $\langle p_n|\enskip n\in N\rangle$ is the sequence of elements of the set of all primes $P$.
   Then it is easy to see that $a_n$ satisfies conditions of the previous theorem
    so  for any $r\in R^+$ there is
   $S\subseteq P$ such that $\sum_{n\in S} -\log(1-1/p_n)=r$,
   i.e. $\prod_{n\in S} (1-1/p_n)=e^{-r}$. As $r$ runs over $R^+$,
   $e^{-r}$ takes all values in $[0,1]$, so for any
   $t\in [0,1]$, there is $S\subseteq P$ such that
   $\prod_{p\in S} (1-1/p)= t$. Thus
    $\lim_{n\in S} \varphi(n)/n= t$, hence
    $A$ is dense in [0,1].
\hfill $\square$
\vspace{3ex}

The next result is useful in studying of dense-expandability.
\begin{thm}
\label{7}
Suppose that $E$ is a sequence of positive real numbers such that $\lim E=0$ and $\sum E=\infty$.
Then for any nonnegative real number $r$
the sequence $r+E$ is dense-expandable.
\end{thm}
\noindent{\bf Proof.} For fixed $r\geqslant0$ we want to prove that $\overline{
\operatorname{X}(r+E)}=\mathbb{R}$. Since $x\in \operatorname{X}(r+E)$
iff $-x\in \operatorname{X}(r+E)$, it is sufficient to prove
that for any $c\geqslant0$ and an arbitrary small $\varepsilon>0$
open interval $(c-\varepsilon,c+\varepsilon)$ and $\operatorname{X}(r+E)$  meet
each other.

The assumed properties of $E$ provide the existence of positive
integers $n$ and $m>n$ such that:
\begin{enumerate}
\item $E(i)<\varepsilon/2$, for all $i\geqslant n$;
\item $\sum\limits_{i=n}^mE(i)\leqslant c<\sum\limits_{i=n}^{m+1}E(i)<c+\varepsilon$.
\end{enumerate}
Let $q=\sum\limits_{i=n}^{m+1}E(i) -c$, $l=m-n+2$ and let $\delta=\displaystyle\frac{q}{2l}$.
Since $\lim E=0$, there is an integer $k>m$ such that $E(i)<\delta$ for all
$i\geqslant k$. Then:
\begin{eqnarray*}
c+\varepsilon&>&\sum_{i=n}^{m+1}(r+E(i))-\sum_{i=k}^{k+l-1}(r+E(i))\\
&=& \sum_{i=n}^{m+1}E(i)-\sum_{i=k}^{k+l-1}E(i)
\geqslant \sum_{i=n}^{m+1}E(i)-\frac{q}{2}
>c.
\end{eqnarray*}
Finally,
$
\sum\limits_{i=n}^{m+1}(r+E(i))-\sum\limits_{i=k}^{k+l-1}(r+E(i))\in
\operatorname{X}(r+E),
$
so $\operatorname{X}(r+E)$ is dense in $\mathbb{R}$.
\hfill $\square$
\vspace{3ex}

An immediate consequence of theorem \ref{7} is the fact that being dense-expandable
is not invariant to asymptotic equivalence. For instance, sequences $E_1=\langle
1\ |\ n\in\mathbb{N}\rangle$ and $E_2=\langle 1+ \frac{1}{n+1}\ |\ n\in\mathbb{N}
\rangle$ are asymptotically equivalent, but $\overline{\operatorname{X}(E_1)}=
\mathbb{Z}$ and $\overline{\operatorname{X}(E_2)}=\mathbb{R}$.
\vspace{2ex}

In general, a cofinite subsequence of a dense-expandable sequence $E$ need not
 be dense expandable. As we have mentioned earlier, the sequence
$$
E=\langle 1,2,2^{-1}, 2^2,2^{-2}, 2^3,2^{-3},\dots\rangle
$$
is dense-expandable, but its cofinite subsequence
$E_1=\langle E(n+2)\ |\ n\in\mathbb{N}\rangle$ is not since
$\operatorname{X}(E_1)\cap (1,2)=\emptyset$.
\vspace{2ex}

The basic strategy in proving that a certain sequence $E$ is dense-expandable
is in choosing countably many pairwise disjoint finite subsets $S_n$ of
$\mathbb{N}$ and appropriate $\varepsilon_{n,i}$s such that the sequence
$$
\langle \sum_{i\in S_n}\varepsilon_{n,i}E(i)\ |\ n\in\mathbb{N}\rangle
$$
satisfies conditions of theorem \ref{7}. As an illustration we will prove that
sequence $E=\langle \ln n\ |\ n>0\rangle$ is dense-expandable. First, note
that the sequence $\langle \ln (1+\frac{1}{2n})\ |\ n>0\rangle$ satisfies the conditions
of theorem \ref{7}, so it is dense-expandable. The sets $S_n=\{2n,2n+1\}$, $n>0$
are pairwise disjoint and
$$
\ln (1+\frac{1}{2n})=\ln(2n+1) - \ln n,
$$
so $\operatorname{X}(\langle \ln (1+\frac{1}{2n})\ |\ n>0\rangle)\subseteq
\operatorname{X}(E)$. Hence $E$ is dense-expandable.
\begin{thm}
The sequence $\langle n^\delta\ |\ n\in\mathbb{N}^+\rangle$ is dense-expandable
if and only if $\delta=-1$ or $\delta>-1$ and $\delta\notin \mathbb{Z}$.
\end{thm}
{\bf Proof} If $\delta<-1$, then $\sum\limits_{n=1}^\infty n^\delta$ converges,
so $\operatorname{X}(\langle n^\delta\ |\ n\in\mathbb{N}^+\rangle)$ is
bounded in $\mathbb{R}$. By theorem \ref{7} sequence $\langle n^\delta\ |\ n\in\mathbb{N}^+\rangle$
is dense-expandable for each $\delta\in [-1,0)$. If $\delta$ is a positive integer,
then $\operatorname{X}(\langle n^\delta\ |\ n\in\mathbb{N}^+\rangle)
\subseteq \mathbb{Z}$. It remains to prove that $\langle n^\delta\ |\ n\in\mathbb{N}^+\rangle$
is dense-expandable for any $\delta\in\mathbb{R}^+\setminus\mathbb{Z}$.
\vspace{2ex}

Fix $\delta\in\mathbb{R}^+\setminus\mathbb{Z}$. Then there is a unique
positive integer $m$ such that $m-1<\delta<m$. Let
$$
E(n)=_{\rm def} \sum_{i=0}^k a_i(n-i)^\delta,\ \ n>k,
$$
where $P_m=\langle a_0,\dots,a_k\rangle$ is the $m$-th $P$-sequence.
Then:
\begin{eqnarray*}
E(n)&=& \sum_{i=0}^k a_i(n-i)^\delta\\
&=& n^\delta\sum_{i=0}^ka_i\left(1- \frac{i}{n}\right)^\delta\\
&=& n^\delta \sum_{i=0}^ka_i\left(\sum_{j=0}^\infty (-1)^j\binom{\delta}{j}i^jn^{-j}\right)\\
&=& n^\delta \sum_{j=0}^\infty (-1)^jn^{-j}\binom{\delta}{j}
\left(\sum_{i=0}^k a_ii^j\right)\\
&=& n^\delta \sum_{j=m}^\infty (-1)^jn^{-j}\binom{\delta}{j}
\left(\sum_{i=0}^k a_ii^j\right)\ \ \ (\mbox{corollary \ref{5}})\\
&=& (-1)^m n^{\delta-m}\binom{\delta}{m}\sum_{i=0}^k a_ii^m\ \ +\ \ o(n^{\delta-m}).
\end{eqnarray*}
Since $\binom{\delta}{m}>0$ and $\operatorname{sgn}\sum\limits_{i=0}^ka_ii^m=
(-1)^m$, sequence $\langle E(n)\ |n>k\rangle$ is ultimately positive.
Now $-1<\delta-m<0$ implies that $\langle E(n)\ |n>k\rangle$ is dense-expandable
(theorem \ref{7}). The same is obviously true for the sequence
$$
E=\langle E(kn)\ |\ n>k\rangle.
$$
Finally, sets $S_n=\{kn-i |\ i\in\{0,\dots,k\}\}$, $n>k$ are pairwise
disjoint and each $E(kn)$ is equal to $\sum\limits_{j\in S_n}\varepsilon_{n,j}
j^\delta$, where $\varepsilon_{n,j}$ are the corresponding coordinates of $P_m$.
Thus  the sequence $\langle n^\delta\ |\ n>0\rangle$ is dense-expandable.
\hfill $\square$
\begin{thm}
Let $\langle p_n\ |\ n\in\mathbb{N}\rangle$ be the sequence of all
prime numbers and let $E=\langle p_n^\delta\ |\ n\in\mathbb{N}\rangle$.
Then:
\begin{enumerate}
\item The Riemann hypothesis implies that $E$ is dense-expandable for
any $0<\delta<1/2$;
\item Hypothesis $\lim\limits_{n\to\infty}(\sqrt{p_{n+1}}-\sqrt{p_n})=0$ implies that
$E$ is dense-expandable for any $0<\delta\leqslant 1/2$.
\end{enumerate}
\end{thm}
{\bf Proof} In order to prove (1), assume the Riemann hypothesis. Then, the following relation
holds for the consecutive primes:
\begin{eqnarray}
p_{n+1}-p_n\ll \sqrt{p_n}\log p_n
\end{eqnarray}
Suppose that $0<\delta<\frac12$. Then
\begin{eqnarray*}
         p_{n+1}^{\delta}-p_n^{\delta}&=&
               p_{n+1}^{\delta}\left[1-\left(1-\frac{p_{n+1}-p_n}{p_{n+1}}\right)^{\delta}\right]  \\
                                      &=& p_{n+1}^{\delta}\left[1-
                                         \left(1-\binom{\delta}{1}\frac{p_{n+1}-p_n}{p_{n+1}}+
                                         o\left(\frac{p_{n+1}-p_n}{p_{n+1}}\right)\right)\right]   \\
                                      &=& \delta\frac{p_{n+1}-p_n}{p_{n+1}^{1-\delta}}+
                                         o\left(\frac{p_{n+1}-p_n}{p_{n+1}^{1-\delta}}\right) \ll
                                         \frac{\sqrt{p_n}\log p_n}{p_{n+1}^{1-\delta}}+o(1)\to 0
\end{eqnarray*}
as $n\to\infty$. Taking $u_n= p_{n+1}^{\delta}-p_n^{\delta}$, $n\in N^+$, we see that for $0<\delta<\frac12$,
$$
\sum_{k\leq n} u_k= p_{n+1}^{\delta}-2^{\delta}\to \infty,\hskip 6mm \text{as}\hskip 3mm n\to \infty,
$$
$u_n>0$ and $\lim\limits_n u_n=0$. Thus, by Theorem \ref{7} $E$ is dense-expandable.
\vspace{2ex}

In order to prove (2), let us assume the hypothesis
$\lim\limits_{n\to\infty} (\sqrt{p_{n+1}}-\sqrt{p_n})=0$. If
$f(\delta)= x^\delta -y^\delta$, $0<\delta$ and $y<x$, then
$f'(\delta)= x^\delta\ln(x) - y^\delta\ln(y)>0$, so $f(\delta)$ is
increasing for  $\delta>0$. Hence, if $0<\delta\leq\frac{1}{2}$,
then $0\leq p_{n+1}^{\delta}-p_n^{\delta}\leq \sqrt{p_{n+1}}-\sqrt{p_n}\to 0$, as
$n\to\infty$, so by an argument as in (1), the assertion follows.
\hfill $ \square$

\section{Sums of like powers}

Finite disjoint subsets $U$ and $V$ of $\mathbb{Z}$ will be called
a Prouhet-Tarry-Escott pair for the given integer $n>1$ if they
have the same cardinality and
\begin{eqnarray}
\label{9}
\sum_{u\in U}u^s=\sum_{v\in V}v^s,\ \ s=0,\dots,n-1,\ \ {\rm and}\
\ \sum_{u\in U}u^n\not=\sum_{v\in V}v^n.
\end{eqnarray}
The sums satisfying the left hand conjunct of (\ref{9}) are also known as sums of like powers.

If $\langle U_1,V_1\rangle,\dots,\langle U_m,V_m\rangle$ are
Prouhet-Tarry-Escott pairs for the given integer $n$ and if sets $U_1,\dots,U_m,V_1,
\dots,V_m$ are pairwise disjoint, then clearly sets
$U=\bigcup\limits_{i=1}^m U_i$ and $V=\bigcup\limits_{i=1}^m$
form another Prouhet-Tary-Escott pair for $n$.
\vspace{2ex}

Now let us describe how one can use the $P$-sequences in order to
generate the Prouhet-Tarry-Escott pairs:
\vspace{2ex}

\noindent $\mathbf 1$ Let $n\geqslant 2$ be an arbitrary integer and let
$P_n=\langle a_0,\dots ,a_k\rangle$ be the $n$-th $P$-sequence. By
lemma \ref{1}  we have that
$$
\sum_{i=0}^ka_i(pi+l)^s=0,\ \ s=0,\dots,n-1, \ l\in\mathbb{Z},\ \
p\in\mathbb{Z}\setminus\{0\}
$$
(observe that $\sum\limits_{i=0}^ka_i(pi+l)^s=p^sF_{n,s}(l/p)$).
Since each $P$-sequence has the same number of $1$s and $-1$s, we
have that sets $U_{p,l}$ and $V_{p,l}$ defined by
$$
U_{p,l}=\{pi+l\ |\ 0\leqslant i\leqslant k\land a_i=-1\}\ {\rm
and}\ V_{p,l}=\{pi+l\ |\ 0\leqslant i\leqslant k\land a_i=1\}
$$
form a Prouhet-Tarry-Escott pair for the given integer $n\geqslant 2$.
\vspace{2ex}

\noindent $\mathbf 2$ Let $P_n=\langle a_0,\dots,a_k\rangle$ be
the $n$-th $P$-sequence ($n\geqslant 2$). We define
the sequence $Q_n=\langle b_0,\dots, b_{k+2}\rangle$ as
follows:
$$
b_i=\left\{
\begin{array}{ccc}
a_i&,&i\in\{0,1,k+1,k+2\}\\
a_i+a_{i+2}&,& \mbox{otherwise}
\end{array}
\right.
.$$
For example, we obtain $Q_3$ from $P_3=\langle1,-1,-1,0,1,1,-1\rangle$
in the following manner:
$$
\begin{array}{rrrrrrrrrrr}
P_3& &1&-1&-1&0&1&1&-1&\\
& & & &1&-1&-1&0&1&1&-1\\ \hline
Q_3& &1&-1&0&-1&0&1&0&1&-1
\end{array}
.$$
An easy induction argument yields that each $b_i\in\{-1,0,1\}$
and that each $Q_n$ has the same number of 1s and -1s. Now for any
non-negative integer $s<n$ we have that
$$
\sum_{i=0}^kb_i(i+1)^s=\sum_{i=0}^ka_i(i+1)^s\ +\ \sum_{i=0}^ka_i(i+3)^s=0,
$$
so $U=\{i+1\ |\ b_i=1\}$ and $V=\{i+1\ |\ b_i=-1\}$ represents a Prouhet-Tarry-Escott pair.
\vspace{2ex}

\noindent $\mathbf 3$
For the $n$-th $P$-sequence $P_n=\langle
a_0,\dots,a_k\rangle$ let
$$
X_n=\{i\in\mathbb{N}^+\ |\ i\leqslant k\land a_i=-1\}\ {\rm and}\
Y_n=\{i\in\mathbb{N}^+\ |\ i\leqslant k\land a_i=1\}.
$$
Clearly, $X_n$ and $Y_n$ are disjoint and $|X_n|=|Y_n|+1$.
Furthermore, using the definition of the notion of a $P$-sequence
one can easily check that
$$
X_{2n+1}\subset X_{2n+2}\ {\rm and}\ Y_{2n+1}\subset Y_{2n+2},
$$
and the sets $U=X_{2n+2}\setminus X_{2n+1}$ and
$V=Y_{2n+2}\setminus Y_{2n+1}$ are disjoint and have the same
cardinality. Bearing in mind the corollary 1, we see that for each
nonnegative integer $s\leqslant 2n$ holds
$$
\sum_{i\in U}i^s=\sum_{i\in X_{2n+2}}i^s-\sum_{i\in X_{2n+1}}i^s=
\sum_{i\in Y_{2n+2}}i^s-\sum_{i\in Y_{2n+1}}i^s=\sum_{i\in V}i^s.
$$
For instance, if $n=4$, then $X_4=\{1,2,6,7,11,12\}$,
$Y_4=\{4,5,8,9,13\}$, $U=\{7,11,12\}$, $V=\{8,9,13\}$ and
$$
7^s+11^s+12^s=8^s+9^s+13^s,\ \ s=1,2.
$$


\begin{thebibliography}{55}
\bibitem{Borwein} P. Borwein, C. Ingalls,
The Prouhet-Tarry-Escott Problem Revisited, {\it Enseign. Math.}
40(1994), 3-27
\bibitem{Borwein2} P. Borwein, P. Lisonek, C. Percival,
Computational Investigations of the Prouhet-Tarry-Escott Problem,
 {\it Math. Comp.}, vol.72, 244(2002), 2063-2070
\bibitem{Dickson} L. E. Dickson, {\it History of the Theory of Numbers},
Carnegie Institute, Washington, 1920, Vol II, Ch. 24
\bibitem{Gloden} A. Gloden, {\it Mehrgradige Glaichungen}, Noordhoff, Groningen, 1944, 103pp.
\bibitem{Guy} R. K. Guy, {\it Unsolved Problems in Number Theory}, Sections A8
({\it Gaps between primes. Twin primes}) and D1 ({\it Sums of like
powers. Euler's conjecture}), second edition, Springer-Verlag, New
York, 1994
\bibitem{Ivic} A. Ivi\' c, {\it An introduction to analytic number theory},
 Izd. knj. Zorana Stojanovi\' ca, Novi Sad, 1996, 389pp.
\bibitem{Narkiewicz} W. Narkiewicz, {\it Classical problems in number theory}, PWN,
Warszawa, 1986, 363pp.
\bibitem{} P. Ribenboim, {\it My Numbers, My Friemds}, Springer-Verlag, New York, 2000, 307pp.
\end{thebibliography}
\end{document}